\documentclass{amsart}
\usepackage{amssymb,latexsym,amsmath}
\usepackage[mathscr]{eucal}
\usepackage{ulem}
\usepackage{verbatim}

\def\sideremark#1{\ifvmode\leavevmode\fi\vadjust{\vbox to0pt{\vss
\hbox to 0pt{\hskip\hsize\hskip1em
\vbox{\hsize2cm\tiny\raggedright\pretolerance10000
\noindent#1\hfill}\hss}\vbox to8pt{\vfil}\vss}}}

\newtheorem*{theorem*}{Theorem}
\newtheorem*{corollary*}{Corollary}
\newtheorem{theorem}{Theorem}[section]
\newtheorem{corollary}[theorem]{Corollary}
\newtheorem{lemma}[theorem]{Lemma}
\newtheorem{proposition}[theorem]{Proposition}
\newtheorem{remark}[theorem]{Remark}
\newtheorem{definition}[theorem]{Definition}
\newtheorem{example}[theorem]{Example}
\newtheorem{question}[theorem]{Question}
\newcommand{\be}{\begin{equation}\label}
\newcommand{\ee}{\end{equation}}
\newcommand{\bq}{\begin{equation*}}
\newcommand{\eq}{\end{equation*}}
\newcommand{\ba}{\begin{align*}}
\newcommand{\ea}{\end{align*}}
\newcommand{\bp}{\begin{proof}}
\newcommand{\ep}{\end{proof}}
\newcommand{\bL}{\begin{lemma}\label}
\newcommand{\eL}{\end{lemma}}
\newcommand{\bP}{\begin{proposition}\label}
\newcommand{\eP}{\end{proposition}}
\newcommand{\bC}{\begin{corollary}\label}
\newcommand{\eC}{\end{corollary}}
\newcommand{\bT}{\begin{theorem}\label}
\newcommand{\eT}{\end{theorem}}
\newcommand{\bTT}{\begin{theorem*}\label}
\newcommand{\eTT}{\end{theorem*}}
\newcommand{\bR}{\begin{remark}\label}
\newcommand{\eR}{\end{remark}}
\newcommand{\bD}{\begin{definition}\label}
\newcommand{\eD}{\end{definition}}

\newcommand{\bE}{\begin{example}\label}
\newcommand{\eE}{\end{example}}
\newcommand{\bQ}{\begin{question}\label}
\newcommand{\eQ}{\end{question}}

\DeclareMathOperator{\diag}{diag}

\begin{document}
\title{A SURVEY ON SUBIDEALS OF OPERATORS AND \\AN INTRODUCTION TO SUBIDEAL-TRACES}
\author{Sasmita Patnaik}
\address{University of Cincinnati\\
          Department of Mathematics\\
          Cincinnati, OH, 45221-0025\\
          USA}
\email{sasmita\_19@yahoo.co.in}
\author{Gary Weiss}
\address{University of Cincinnati\\
          Department of Mathematics\\
          Cincinnati, OH, 45221-0025\\
          USA}
\email{gary.weiss@math.uc.edu}
\date{\today}

\dedicatory{\textit{Dedicated to the memory of Mihaly Bakonyi}}
\maketitle

\begin{abstract}  
Operator ideals in $B(H)$ are well understood and exploited but ideals inside them have only recently been studied starting with the 1983 seminal work of Fong and Radjavi and continuing with two recent articles by the authors of this survey.
This article surveys this study embodied in these three articles.
A subideal is a two-sided ideal of $J$ (for specificity also called a $J$-ideal) for $J$ an arbitrary ideal of $B(H)$.
In this terminology we alternatively call $J$ a $B(H)$-ideal. 

This surveys \cite{FR83}, \cite{PW11} and \cite{PW12} in which we developed a complete characterization of all $J$-ideals generated by sets of cardinality strictly less than the cardinality of the continuum. 
So a central theme is the impact of generating sets for subideals on their algebraic structure.
This characterization includes in particular finitely and countably generated $J$-ideals.  
It was obtained by first generalizing to arbitrary principal $J$-ideals the 1983 work of Fong-Radjavi who determined which principal $K(H)$-ideals are also $B(H)$-ideals.
A key property in our investigation turned out to be $J$-softness of a $B(H)$-ideal $I$ inside $J$, that is, $IJ = I$, a generalization of a recent notion of $K(H)$-softness of  $B(H)$-ideals introduced by Kaftal-Weiss and earlier exploited for Banach spaces by Mityagin and Pietsch.
This study of subideals and the study of elementary operators with coefficient constraints are closely related.
Here we also introduce and study a notion of subideal-traces where classical traces (unitarily invariant linear functionals) need not make sense for subideals that are not $B(H)$-ideals.  
\end{abstract}

%%%%%%%%%%%%%%%%%%%%%%%%%%%%%%%%%%%%%%%%%%%%%%%%%%%%%%%%%%%%%%%%%%%%%%%%%%%%%%%%%%%%%%%%%%%%%%%%%%%%%%%%%%%%%%%%%%%%%%%%%%%%%%%%%%%%%%%%%

\section{Introduction}  

For general rings, an ideal (all ideals herein are two-sided ideals) is a commutative additive subgroup of a ring that is closed under left and right multiplication by elements of the ring. 
Herein $H$ denotes a separable infinite-dimensional complex Hilbert space and $B(H)$ denotes the $C^*$-algebra of all bounded linear operators on $H$.
Ideals of $B(H)$ (henceforth alternatively called $B(H)$-ideals) have become ubiquitous throughout operator theory since their celebrated characterization by Calkin and Schatten \cite{C41}, \cite{S}, in terms of ``characteristic sets" of singular number sequences $s(T)$ of the operators $T$ in the ideal.
This characterization has had and continues to have substantial impact in operator theory.
As commutative objects in analysis, characteristic sets make more accessible the subtler properties of $B(H)$-ideals, particularly illuminating and expanding the knowledge of some of their noncommutative features.
Some well-known $B(H)$-ideals are the ideal of compact operators $K(H)$, the finite rank operators $F(H)$, principal ideals $(S)$ (i.e., singly generated $B(H)$-ideals), Banach ideals, the Hilbert-Schmidt class $C_2$, the trace class $C_1$, Orlicz ideals, Marcinkiewicz ideals and Lorentz ideals, to name a few.
Definitions and properties of these ideals among others may be found in \cite{DFWW}.

A subideal of operators is an ideal of $J$ (for specificity called a $J$-ideal) for $J$ an arbitrary $B(H)$-ideal.
``Subideal" is a name coined by Gary Weiss motivated from the 1983 seminal work of Fong-Radjavi and by the new perspectives on operator ideals from work of Dykema, Figiel, Weiss and Wodzicki \cite{DFWW}. 
It is clear that every $B(H)$-ideal is a subideal, but the converse is less clear, i.e., whether or not every subideal is also a $B(H)$-ideal. 
Fong-Radjavi constructed the first example of a principal $K(H)$-ideal that is not a $B(H)$-ideal (Example \ref{E:1.1}).
This shows that the class of subideals is strictly larger than the class of $B(H)$-ideals.

The main and most general results in this survey are Theorem \ref{T:2.1} and Theorem \ref{T:3} (Structure theorem for subideals $(\mathscr S)_J$ for $|\mathscr S| < \mathfrak c$) in which we characterize in terms of a new notion called softness, when a subideal generated by strictly less than $\mathfrak c$ elements is also a $B(H)$-ideal; and then we characterize its algebraic structure.
Section 4 compares $B(H)$-ideals to subideals via some of their differences and similarities.
And Section 5 is new research that begins the investigation of subideal-traces, an attempt and an analog to traces on $B(H)$-ideals which are themselves ubiquitous in operator theory.

%%%%%%%%%%%%%%%%%%%%%%%%%%%%%%%%%%%%%%%%%%%%%%%%%%%%%%%%%%%%%%%%%%%%%%%%%%%%%%%%%%%%%%%%%%%%%%%%%%%%%%%%%%%%%%%%%%%%%%%%%%%%%%%%%%%%%%%%%%%%

\section{Preliminaries}

Every $B(H)$-ideal $J$ is linear because for each $\alpha \in \mathbb C$, $\alpha {\bf1} \in B(H)$, so then for each $A \in J$, $\alpha A=(\alpha {\bf1})A \in J$.
But surprisingly a subideal (i.e., a $J$-ideal) may not be linear (Section 4-Example \ref{E:3}, see also \cite[Example 3.5]{PW11}). 
In  \textit{Subideals of Operators} \cite{PW11} we found three types of principal and finitely generated subideals: linear, real-linear and nonlinear classical subideals. 
Such types also carry over to non-principal $J$-ideals.
The linear $K(H)$-ideals, being the traditional linear ones, were studied in 1983 by Fong-Radjavi \cite{FR83}. 
They found principal linear $K(H)$-ideals that are not $B(H)$-ideals.
Herein we take all $J$-ideals to be linear, but as proved in \cite{PW11}, we expect here also that most of the results and methods apply to the two other types of subideals (real-linear and nonlinear classical).

Noting the obvious fact that intersections of ideals in any ring are themselves ideals, we begin with the following definition.  

\bD{D:0}
\item{(i)} The principal $B(H)$-ideal generated by the single operator $S$ is defined by
\begin{center}
$\left(S\right)$ := $\bigcap \{I \mid I$ is a $B(H)$-ideal containing S$\}$
\end{center}
\item{(ii)} The principal $J$-ideal generated by $S$ is defined by
\begin{center}
$\left(S\right)_J$ := $\bigcap \{I \mid  I$ is a $J$-ideal containing $S\}$
\end{center}

\item{(iii)} As above for principal $J$-ideals, likewise for an arbitrary subset $\mathscr S \subset J$, $(\mathscr S)$ and  $(\mathscr S)_J$ denote respectively, 

via intersections, the smallest $B(H)$-ideal and the smallest $J$-ideal generated by the set $\mathscr S$.
\eD

\bD{D:10} \quad \\
For $B(H)$-ideals $I, J$, ideal $I$ is called ``$J$-soft" if $IJ = I$. (Clearly this applies only when $I\subset J$.)  \\
Equivalently in the language of s-numbers (see Remark \ref{R:11}(i),(ii),(v) below): \\
For every $A \in I$, $s_n(A) = O(s_n(B)s_n(C))$ for some $B \in I, C \in J$. \\
($s(A) := \left<s_n(A)\right>$ is the singular number (s-number) sequence of operator $A$, counting multiplicities of course.)
\eD

\bR{R:11} \emph{Standard facts and tools for operator ideals.}\\

(i) If $I,J$ are $B(H)$-ideals, then the ideal product $IJ$, which is both associative and commutative, is the $B(H)$-ideal alternatively defined via its characteristic set is given by $\Sigma(IJ) = \{\xi \in c_o^* \mid \xi \leq \eta\rho$ for some $\eta \in \Sigma(I)$ and $\rho \in \Sigma(J)\}$ \cite[Sections 2.8, 4.3]{DFWW} (see also \cite[Section 4]{KW07}).(See also Historical Background below-first paragraph.)\\

(ii) If $I$ and $J$ are $B(H)$-ideals for which $A \in IJ$, then $A = XY$ for some $X \in I, Y \in J$ \cite[Lemma 6.3]{DFWW}.\\

(iii) For $T \in B(H)$, $ A \in (T)$ if and only if $s(A) = \textrm{O}(D_m(s(T)))$ for some $m \in \mathbb{N}$.
Moreover, for $I$ a $B(H)$-ideal, $A \in I$ if and only if $A^* \in I$ if and only if $|A| \in I$ (via the polar decomposition $A = U|A|$ and $U^*A = |A| := (A^*A)^{1/2}$) if and only if $\diag s(A) \in I$. $D_m \xi$ is the $m$-fold sequence ampliation recalled just below in Historical Background.\\

(iv) The lattice of $B(H)$-ideals forms a commutative semiring with multiplicative identity $B(H)$.
That is, this lattice is commutative and associative under ideal addition and multiplication (see \cite[Section 2.8]{DFWW}) and it is distributive.
Distributivity with multiplier $K(H)$ is stated without proof in \cite[Lemma 5.6-preceding comments]{KW07}.
 
One importance of principal ideals in a general ring is that they are building blocks for all ideals $I$ that contain them in that: 
$I \quad = \displaystyle{\bigcup_{r_1,\dots,r_n \in I, \,n \in \mathbb N}} (r_1)+ \cdots +(r_n).$
Note also $(r) = r + Rr + rR + \sum_{\text{finite sum}} RrR$, and when $R$ is unital, then $(r) = Rr + rR + \sum_{\text{finite sum}}RrR$.

Finally when $R = J$ is a  $B(H)$-ideal, $(r)$ collapses to $(r) = r + Rr + rR + RrR$ because as proved in \cite[Lemma 6.3]{DFWW}, $\sum_{\text{finite sum}} RrR = RrR$. 
\\

(v) When $T = \displaystyle{\sum_{i=1}^{n}}A_{i}TB_{i}$ with each $A_i$ or $B_i \in J$, the important s-number relation holds:\\
\noindent \qquad $s(T) = \textrm{O}(D_{m}(s(T))s(C))$ for some $C \in J$ (since then $T \in (T)J$, see \cite[Section 1, p. 6]{KW07} and Remark \ref{R:11}(i)).\\
\eR

\textbf{\textit{Historical Background}}\\

Calkin-Schatten completely characterized $B(H)$-ideals via the lattice preserving isomorphism between $B(H)$-ideals and characteristic sets $\Sigma \subseteq c_{0}^*$ where %$I \rightarrow \Sigma(I)$ induced by $I \owns X \rightarrow s(X) \in \Sigma(I)$.
$c_{0}^{*}$ denotes the cone of nonnegative sequences decreasing to zero; characteristic sets $\Sigma$ are those subsets of $c_{0}^{*}$ that are additive, hereditary (solid) and ampliation invariant (invariant under each $m$-fold ampliation $D_m\xi :=\,<\xi_{1},\cdots, \xi_{1}, \xi_{2}, \cdots, \xi_{2}, \cdots>$ with each entry $\xi_{i}$ repeated $m$ times);
the characteristic set $\Sigma(I) := \{\eta \in   c_{0}^* \mid \diag \eta \in I\}$,
so for example $\Sigma (K(H)) = c_{0}^*$.

In 1983 Fong-Radjavi \cite{FR83} investigated principal $K(H)$-ideals.
They found principal $K(H)$-ideals that are not $B(H)$-ideals (Example \ref{E:1.1} below)
by determining necessary and sufficient conditions for a  principal $K(H)$-ideal to be also a $B(H)$-ideal
\cite[Theorem 2]{FR83}.
And in doing so, at least for the authors of this paper, they initiated the study of subideals.
The main results of Fong-Radjavi are summarized in the following theorem. 
 
\begin{theorem*} \cite[Theorems 1-2]{FR83} For $T$ a compact operator of infinite rank, 
$P := (T^*T)^{\frac{1}{2}}$ and $\mathcal{I}$ the ideal in $K(H)$ generated by $T$ 
and $\mathcal P$ the ideal of $K(H)$ generated by $P$, the following are equivalent.

(i) $\mathcal{I}$ is an ideal in $B(H)$.

(ii) $\mathcal{P}$ is an ideal in $B(H)$.

(iii) $\mathcal{P}$ is a Lie ideal in $B(H)$.

(iv) $T =   A_1TB_1 + \dots + A_kTB_k$ for some $k$ and some $A_i \in K(H)$, $B_i \in B(H)$.

(v) $T =   A_1TB_1 + \dots + A_kTB_k$ for some $k$ and some $A_i, B_i \in K(H)$.

(vi) For some integer $k > 1$, $s_{nk}(P)= o(s_n(P))$ as $n \rightarrow \infty$.
\end{theorem*}
Fong-Radjavi proved this via the positive case employing the Lie ideal condition (iii),
but our approach below avoids the need for considering separately the positive case and any Lie ideal considerations.
Notably also, conditions (iv)-(v) above indicates the relevance of elementary operators with coefficient constraints. 
 
\bE{E:1.1}
Condition (vi) of  the above theorem shows that if the singular number sequence of the operator $P$ is given by $s(P)= \left<\frac{1}{2^n}\right>$,  then the principal $K(H)$-ideal generated by $P$ is a $B(H)$-ideal.
But if $s(P)= \left<\frac{1}{n}\right>$, then the principal $K(H)$-ideal generated by $P$ is not a $B(H)$-ideal.\\
Indeed, $\frac{\frac{1}{2^{nk}}}{\frac{1}{2^n}} = \frac{1}{2^{n(k-1)}} \rightarrow 0$ but $\frac{\frac{1}{nk}}{\frac{1}{n}}= \frac{1}{k} \nrightarrow 0$ as $n \rightarrow \infty$.
\eE

%%%%%%%%%%%%%%%%%%%%%%%%%%%%%%%%%%%%%%%%%%%%%%%%%%%%%%%%%%%%%%%%%%%%%%%%%%%%%%%%%%%%%%%%%%%%%%%%%%%%%%%%%%%%%%%%%%%%%%%%%%%%%%%%%%%%%%%%%

\section{Subideals of Operators}
\noindent Motivated by the Calkin-Schatten characterization and the seminal work of Fong-Radjavi, a natural question to ask is:
\begin{center}
\textit{What can be said about subideals, i.e., is it possible to characterize them in some way?}
\end{center}
\vspace{.2cm}

A conventional approach to attack the characterization problem for $J$-ideals is to begin at the elementary level as did Fong-Radjavi, albeit they did not consider characterizations except implicitly for principal $K(H)$-ideals in one of their proofs.
So we first investigate principal $J$-ideals, 
then finitely generated $J$-ideals and then $J$-ideals  $\mathcal I = (\mathscr S)_J$ generated by sets $\mathscr S$ of higher cardinalities including the countable case.  
We fully generalize Fong-Radjavi's  \cite[Theorem 2]{FR83} from principal $K(H)$-ideals to arbitrary principal $J$-ideals and then to finitely generated $J$-ideals. 
The reason to consider the finitely generated case separate from the principal case is that, unlike $B(H)$-ideals where every finitely generated $B(H)$-ideal is always a principal $B(H)$-ideal, a finitely generated $J$-ideal need not be a principal $J$-ideal (see Section 4, Example \ref{E:4} for the case $J = K(H)$).
Consequently, we characterize all $J$-ideals generated by sets of cardinality strictly less than the cardinality of the continuum, including finitely and countably generated $J$-ideals.
A key property in this characterization turned out to be $J$-softness of a $B(H)$-ideal $I$ inside $J$, that is, $IJ = I$ (Definition \ref{D:10}) a generalization of a recent notion of $K(H)$-softness of $B(H)$-ideals introduced by  Kaftal-Weiss \cite{KW07} and earlier exploited for Banach and Hilbert spaces by Mityagin and Pietsch.\\

\newpage
We first begin with the following algebraic description of the principal $J$-ideal generated by $S \in J$ (see Remark \ref{R:11}(iv)).

\bP{P:11}
 For $S \in J$, an algebraic description of principal $J$-ideal $(S)_J$ is given by 
\begin{equation*}
(S)_J = \{\alpha S + AS + SB + \displaystyle{\sum_{i=1}^{m}}A_iSB_i \mid  A,\,B,\,A_i, \,\,B_i \in J,\,\alpha \in \mathbb{C},\,m \in \mathbb{N}\}
\end{equation*}
That is, $(S)_J = \mathbb CS + JS+ SJ+ J(S)J$.
\eP

The following theorem generalizes Fong-Radjavi's result from principal $K(H)$-ideals to principal $J$-ideals by determining necessary and sufficient conditions for a principal $J$-ideal to be also a $B(H)$-ideal. Here is where $J$-softness first played a prominent role.
 \vspace{.2cm}

For compact operators $S, T$, the product $s(S)s(T)$ denotes the pointwise product of their s-number sequences.

\bT{T:3.2.2}
For $S \in J$ and $(S)_J$, the principal $J$-ideal generated by $S$, the following are equivalent.
\item (i) $(S)_J$ is a $B(H)$-ideal.

\item (ii) The principal $B(H)$-ideal
$(S)$ is $J$-soft, i.e., $(S)$ = $J(S)$ (equivalently, $(S)$ = $(S)J$).

\item (iii) $S = AS + SB + \displaystyle{\sum_{i=1}^{m}}A_{i}SB_{i}$ for some $A,\,B,\,A_{i},\, B_{i} \in J,~ m \in \mathbb{N}$.

\item(iv) $s(S) = \textrm{O}(D_{k}(s(S))s(T))$ for some $T \in J$ and $ k \in \mathbb{N}$.\\
\eT

\bp[Proof of $(i) \Rightarrow (ii)$ only]
This is the main part of the proof so we provide here an outline.
For every unitary map $\phi: H \rightarrow H \oplus H$, $S \mapsto \phi S \phi^{-1}$ preserves s-number sequences and hence also ideals via Calkin-Schatten's representation.
Since $(S)_{J}$ is a $B(H)$-ideal containing $S$, $\phi^{-1}(S \oplus 0)\phi,\, \phi^{-1}(0 \oplus S)\phi \in (S)_{J}$ since they possess the same s-numbers as $S$.
Then by Proposition \ref{P:11} for principal $J$-ideal $(S)_J$,
$$\phi^{-1}(S \oplus 0)\phi = \alpha S+ X$$  $$\phi^{-1}(0 \oplus S)\phi = \beta S+ Y$$ for $X, Y \in JS + SJ + J(S)J$.

If $\alpha$ = 0 or $\beta = 0$, then $\phi^{-1}(S \oplus 0)\phi$ or $\phi^{-1}(0 \oplus S)\phi \in J(S)$.
Then, in either case, $S \in J(S)$, hence $(S) \subseteq J(S)$ and since the other inclusion is automatic, one has $(S) = J(S)$.
If $\alpha, \beta \ne 0$, multiplying the first equation by $-\beta$ and the second equation by $\alpha$ and adding obtains
$\phi^{-1}(-\beta S \oplus \alpha S)\phi = -\beta X + \alpha Y \in J(S)$.
Multiplying $-\beta S \oplus \alpha S$ in $B(H \oplus H)$ by a suitable diagonal projection one obtains $\phi^{-1}(S \oplus 0)\phi \in J(S)$.
Hence, also $S \in J(S)$, again equivalent to (ii). 
\ep

\bR{R:1}
Using basic linear algebra techniques, we extended Theorem \ref{T:3.2.2} from principal $J$-ideals to finitely generated $J$-ideals by solving a large system of linear equations which we then project into a finite dimensional quotient space \cite[Theorem 4.5]{PW12}.
\eR

The techniques for finitely generated subideals do not work for countably generated subideals because then the latter case involves an intractable infinite system of equations, so a more sophisticated approach was needed.
Based on the Hamel dimension of a related quotient space, a necessary and sufficient softness condition is found for subideals to also be $B(H)$-ideals among those subideals with generating sets of cardinality strictly less than $\mathfrak c$, so includes all countably generated subideals (Theorem \ref{T:2.1}, see also \cite[Theorem 4.1]{PW12}). 
We then use this condition to characterize these subideals (Theorem \ref{T:3}, see also \cite[Theorem 4.4]{PW12}).
To investigate this in \cite{PW12}, we began with the following proposition.

\bP{P:0.1}\cite[Proposition 3.1]{PW12}
For the $J$-ideal $(\mathscr S)_J$ generated by a set $\mathscr S$ and defining \\$(\mathscr S)_J^0 :=\text{span} \{\mathscr SJ + J\mathscr S\} + J(\mathscr S)J$, the Hamel dimension of the quotient space $(\mathscr S)_J/(\mathscr S)_J^0$ is at most the cardinality of the generating set $\mathscr S$.
\eP

The main softness theorem in \cite{PW12} characterizing when a $J$-ideal is also a $B(H)$-ideal is:

\begin{theorem}\cite[Theorem 4.1]{PW12}\label{T:2.1}
A $J$-ideal $(\mathscr S)_J$ generated by a set $\mathscr S$ of cardinality strictly less than $\mathfrak c$ is a $B(H)$-ideal if and only if  the $B(H)$-ideal $(\mathscr S)$ is $J$-soft.
\end{theorem}

\bp[Sketch of proof]
Here we sketch only the proof of the first implication, that is, that $(\mathscr S)_{J}$ is a $B(H)$-ideal implies $(\mathscr S)$ is $J$-soft. The reverse implication is somewhat routine. 
The algebraic structure of $(\mathscr S)_{ J}$ is given by $(\mathscr S)_{ J} = \text{span}\,\{\mathscr S\} + (\mathscr S)_{ J}^0 $ and so the quotient space $(\mathscr S)_{ J}/(\mathscr S)_{ J}^0 = \text{span}\, \{[S_{\alpha}]\}$ where $S_{\alpha}$ ranges over $\mathscr S$.
Hence the Hamel dimension of $(\mathscr S)_{ J}/(\mathscr S)_{  J}^0$ is strictly less than $\mathfrak c$.
Moreover, since $(\mathscr S)_{  J}$ is also a  $B(H)$-ideal, by minimality $(\mathscr S)_{  J}= (\mathscr S)$. 

The assumption that $(\mathscr S)  J\subsetneq (\mathscr S)$ provides an operator in their difference which we use to construct an imbedding of  $\ell^p$ into $(\mathscr S)_{ J}/(\mathscr S)_{  J}^0$. 
But the Hamel dimension of $\ell^p$ is $\mathfrak c$ \cite[Lemma 3.4]{HN00} and the Hamel dimension of $(\mathscr S)_{ J}/(\mathscr S)_{ J}^0$ is strictly less than $\mathfrak c$, a contradiction.
Therefore, the condition $(\mathscr S)_{ J}$ is a $B(H)$-ideal implies that $(\mathscr S)J = (\mathscr S) $, that is, $(\mathscr S)$ is $J$-soft.
\ep

\bR{R:00} 
Theorem \ref{T:2.1} on the equivalence of a $J$-ideal $(\mathscr S)_J$ being a $B(H)$-ideal and $(\mathscr S)$, the $B(H)$-ideal it generates, being $J$-soft motivates the question on whether this is always true independent of its various classes of generators. 
The answer is no from the following example.
And Theorem \ref{T:2.1} yields new information about the possible cardinality of any  class of its generators.

The $K(H)$-ideal $(\diag \left<\frac{1}{n}\right>)$  is also a principal $B(H)$-ideal but is not $K(H)$-soft \cite[Section 4, Example 4.5]{PW12}.   
Thus  $\mathcal I$ being a $B(H)$-ideal is not equivalent to $J$-softness of the $B(H)$-ideal $(\mathcal I)$, for $\mathcal I$ a $J$-ideal and $(\mathcal I)$ the $B(H)$-ideal generated by $\mathcal I$.
Moreover, by Theorem \ref{T:2.1}, $(\diag \left<\frac{1}{n}\right>)$ which is also a $K(H)$-ideal, cannot be generated in $K(H)$ by less than $\mathfrak c$ generators.
\eR

As a consequence of Theorem \ref{T:2.1} we obtain a characterization of all $J$-ideals generated by sets of cardinality strictly less than the cardinality of the continuum. 
These are the countably generated $J$-ideals when assuming the continuum hypothesis, and otherwise these include more $J$-ideals than the countably generated ones.

\bT{T:3}(Structure theorem for  $(\mathscr S)_J$ when $|\mathscr S| < \mathfrak c$)\\
The algebraic structure of the $J$-ideal $(\mathscr S)_J$ generated by a set $\mathscr S$ of cardinality strictly less than $\mathfrak c$ is given by
\begin{center}
$(\mathscr S)_J = \text{span}\{\mathscr S + J\mathscr S + \mathscr SJ\} + J(\mathscr S)J,$
\end{center}
$J(\mathscr S)J $ is a $B(H)$-ideal, $\text{span}\{J\mathscr S + \mathscr SJ\} + J(\mathscr S)J$ is a $J$-ideal, and $$J(\mathscr S)J \subset \text{span}\{J\mathscr S + \mathscr SJ\} + J(\mathscr S)J \subset (\mathscr S)_J$$
This inclusion collapse to $J(\mathscr S)J = (\mathscr S)_J$ if and only if  $(\mathscr S)$ is $J$-soft (i.e., $(\mathscr S)J = (\mathscr S)$).
\eT

%%%%%%%%%%%%%%%%%%%%%%%%%%%%%%%%%%%%%%%%%%%%%%%%%%%%%%%%%%%%%%%%%%%%%%%%%%%%%%%%%%%%%%%%%%%%%%%%%%%%%%%%%%%%%%%%%%%%%%%%%%%%%%%%%%%%%%%%%%%

\section{Comparison of Subideals to $B(H)$-ideals}

As mentioned in  Preliminaries Section 2, a subideal may not be linear. 
This led the authors of this paper to introduce three kinds of $J$-ideals, namely, linear, real-linear and classical  $J$-ideals (\cite[Definition 2.1]{PW11})(the latter two are nonlinear). 
The term ``classical" is meant in the sense of abstract rings, for instance, ideals where scalar multiplication may not make sense.
The classical principal $J$-ideal generated by $S$ is defined by
$\left<S\right>_J$ := $\bigcap \{I \mid  I$ is a classical $J$-ideal containing $S\}$. 
From Remark \ref{R:11}(iv) one deduces that 
\begin{center}
$\left<S\right>_J = \{ nS + AS + SB + \displaystyle{\sum_{i=1}^{m}}A_iSB_i \mid  A,\,B,\,A_i, \,\,B_i \in J,\, n \in \mathbb{Z},\,m \in \mathbb{N}\}$. 
\end{center}

\bE{E:3}\textit{A concrete nonlinear principal ideal is: $\left<\diag\left<\frac{1}{n}\right>\right>_{K(H)}$}\\
Indeed, if it were linear, then the principal $B(H)$-ideal $(\diag\left<\frac{1}{n}\right>)$ would be $K(H)$-soft, which is not the case. (Combine Example \ref{E:1.1} and Theorem \ref{T:3.2.2}.)
\eE

The explicit description of the principal $J$-ideal generated by $S$ given in Proposition \ref{P:11} implies that
every principal $J$-ideal contains $J(S)J$.
It is well-known that every proper $B(H)$-ideal contains $F(H)$, the $B(H)$-ideal of all finite rank operators \cite[Chapter III, Section 1, Theorem 1.1]{GK69}.
So, every  nonzero principal $J$-ideal contains $F(H)$ (since $S \neq 0$ implies $(S)_J \supset J(S)J \neq \{0\}$) and hence also every nonzero $J$-ideal.
The intersection of all $B(H)$-ideals properly containing $F(H)$ is precisely $F(H)$ \cite[Corollary 3.8(ii)]{KW11},
and since every $B(H)$-ideal is a $J$-ideal, it is clear then that the intersection of all $J$-ideals properly containing $F(H)$ is also precisely $F(H)$.\\

Some striking differences between $J$-ideals and $B(H)$-ideals are described next for the case $J=K(H)$ in Examples \ref{E:4}-\ref{E:0}.
Every finitely generated $B(H)$-ideal is always a principal $B(H)$-ideal because, as is straightforward to see,
the $B(H)$-ideal generated by $\mathscr S = \{S_1, \dots, S_n\} \subset B(H)$, namely $(\mathscr S)$,
is precisely the principal ideal $(|S_1| + \dots + |S_n|)$ where $|S| := (S^*S)^{1/2}$.
But finitely generated $J$-ideals (classical, linear or real-linear) may not be principal as seen in the following example.

\bE{E:4} \textbf{(A doubly generated $J$-ideal of any of the three types that is not principal)}\\
For $J=K(H)$, $S_1 = \diag \, (1, 0, \frac{1}{2}, 0, \frac{1}{3}, \cdots)$ and $S_2 = \diag \, (0, 1, 0, \frac{1}{2}, 0, \frac{1}{3}, \cdots)$, 
$(\{S_1,S_2\})_{K(H)}$ is not a principal linear $K(H)$-ideal, and likewise for the classical and real-linear cases $\left<\{S_{1},S_{2}\}\right>_J$ and $(\{S_{1},S_{2}\})_J^{\mathbb R}$ \cite[Section 4, Example 4.1]{PW11}.
\eE

For $T \in B(H), (T) = (|T|)$, but this need not be true for principal linear $K(H)$-ideals (Example \ref{E:(T) ne (|T|)}).
Moreover, all $B(H)$-ideals are selfadjoint, but this is not necessarily true for principal linear $K(H)$-ideals (Example \ref{E: J ne J^*}) and unlike $B(H)$-ideals, $K(H)$-ideals need not necessarily commute under ideal product (Example \ref{E:0}).

\bE{E:(T) ne (|T|)}
If $J = K(H)$ and operator $T = \diag \left<\frac{(i)^n}{n}\right>$, then $(T)_{K(H)} \neq (|T|)_{K(H)}$.
In fact, $(|T|)_{K(H)}\nsubseteq (T)_{K(H)}$ and $(T)_{K(H)}\nsubseteq (|T|)_{K(H)}$ \cite[Section 5, Example 5.1]{PW11}.
\eE

\bE{E: J ne J^*}(Example of a $K(H)$-ideal that is not closed under the adjoint operation)
$T^{*}\notin (T)_{K(H)}$ where $T = \diag \left<\frac{(i)^n}{n}\right>$ \cite[Section 5, Example 5.2]{PW11}.
\eE
 
\bE{E:0}(Example of $K(H)$-ideals that do not commute)
Consider $J = K(H)$ and with respect to the standard basis
take $S$ to be the diagonal matrix
$S$ := $\text{diag}\,(1, 0, 1/2, 0, 1/3, 0, ...)$
and $T$ to be the weighted shift with this same weight sequence.
Then $(S)_{K(H)}(T)_{K(H)} \neq (T)_{K(H)}(S)_{K(H)}$  \cite[Section 5, Example 5.4]{PW12}.
\eE

%%%%%%%%%%%%%%%%%%%%%%%%%%%%%%%%%%%%%%%%%%%%%%%%%%%%%%%%%%%%%%%%%%%%%%%%%%%%%%%%%%%%%%%%%%%%%%%%%%%%%%%%%%%%%%%%%%%%%%%%%%%%%%%%%%%%%%%%%%

\section{Subideal-Traces}
Subideals $\mathcal I$ that are not $B(H)$-ideals need not be invariant under unitary equivalence, i.e., $U\mathcal IU^* \nsubseteq \mathcal I$ for some unitary operator $U$ (Examples \ref{E:1}-\ref{E:2} below).
Therefore, the definition of trace on a $B(H)$-ideal, that is, a unitarily invariant linear functional, need not make sense on a subideal. Motivated by our work in \cite{BPW12} on unitary operators of the form $U = {\bf 1} + A$ for $A \in K(H)$ we observe that subideals $\mathcal I$ are invariant under these unitaries (i.e., $U\mathcal IU^* \subset \mathcal I$). 
This led the authors of this paper to introduce the notion of a subideal-trace as defined below in Definition \ref{D:1} (see also Remark \ref{R:111}).

\bE{E:1}(Example of a $K(H)$-ideal that is \underline{not} invariant under unitary equivalence)\\
For $J = K(H)$ and a unitary map $\phi: H \rightarrow H \oplus H$, consider $S = \phi^{-1}(D\oplus0)\phi$ for $D = \text{diag} \left<\frac{1}{n}\right>$. Then $(S)_{K(H)}$ the principal $K(H)$-ideal generated by $S$ is not invariant under unitary equivalence. 
We prove this by constructing one unitary operator $U$ for which $USU^* \notin (S)_{K(H)}$.
Indeed, assume $(S)_{K(H)}$ is invariant under unitary equivalence. 
We then have the following contradiction.
 Since \begin{equation*} 
\phi^{-1}\begin{pmatrix}
0  & {\bf1}  \\
{\bf1} &0      
\end{pmatrix}\phi \,\,\text{is a unitary operator in}\,\, B(H),
\end{equation*}
it follows that \begin{equation*}
\phi^{-1}\begin{pmatrix}
0  & {\bf1}  \\
{\bf1} &0      
\end{pmatrix}\phi\,\, S\,\, \phi^{-1}\begin{pmatrix}
0  & {\bf1}  \\
{\bf1} &0      
\end{pmatrix}\phi = \phi^{-1}\begin{pmatrix}
0  & 0  \\
 0& D      
\end{pmatrix}\phi \in (S)_{K(H)} 
\end{equation*}  
 
Using the algebraic structure of $(S)_{K(H)}$ (Proposition \ref{P:11}) one obtains,
\begin{equation*}
\phi^{-1}\begin{pmatrix}
0  & 0  \\
 0& D      
\end{pmatrix}\phi = \alpha S + X,
\end{equation*}
where $ X \in K(H)S + SK(H) + K(H)(S)K(H) \subset (\diag\left<\frac{1}{n}\right>)K(H)$ (since $s(S) = s(D), (S)= (\diag\,\left<\frac{1}{n}\right>)$).
That is, \begin{equation*}
\phi^{-1}\begin{pmatrix}
-\alpha D  & 0  \\
 0& D      
\end{pmatrix}\phi \in (\diag\left<\frac{1}{n}\right>)K(H).
\end{equation*}
This implies that $D \in (\diag\left<\frac{1}{n}\right>)K(H)$, a contradiction to the non-softness of $(\diag\left<\frac{1}{n}\right>)$  \cite[Example 3.3]{PW11}.
Therefore, $(S)_{K(H)}$ is not invariant under unitary equivalence.
\eE

\bE{E:2}(Example of a $K(H)$-ideal that \underline{is} invariant under unitary equivalence) J. Varga  \cite{V} constructed a concrete example of a $K(H)$-ideal generated by the unitary orbit of a positive compact operator that is not a $B(H)$-ideal, namely, $(\mathcal U(A))_{K(H)}$ where $0\leq A \in K(H)$ and $\mathcal U(A) = \{ UAU^* \,|\, U^* = U^{-1}\}$. 
Using Remark \ref{R:11}(iv) for an ideal written as the union of finite sums of its principal ideals, and Proposition \ref{P:11} giving the algebraic structure of the principal $K(H)$-ideal $(UAU^*)_{K(H)}$ generated by $UAU^*$: 
for each $T \in (UAU^*)_{K(H)}$ and $V$ a unitary operator in $B(H)$, from Proposition \ref{P:11} one has
\begin{align*}
VTV^* &= V(\alpha UAU^* + BUAU^* + UAU^*C + A'XB')V^* \qquad \qquad (\text{where}~ B, C, A', B' \in K(H), X \in (UAU^*)) \\ 
&= \alpha VUAU^*V^* + VBV^*VUAU^*V^* + VUAU^*V^*VCV^* + VA'V^*VXV^*VB'V^* \\
&\hspace{8.6cm} (VXV^* \in (VUAU^*V^*) ~since~ X \in (UAU^*))\\
&\in (VUAU^*V^*)_{K(H)} \subset (\mathcal U(A))_{K(H)}\qquad \qquad (since~ VU ~is~ unitary)
\end{align*}
Therefore the $K(H)$-ideal $(\mathcal U(A))_{K(H)}$ is invariant under unitary equivalence.
\eE

Denote by $\mathcal U(H)$ the full group of unitary operators in $B(H)$. 
Recall the essential feature of traces: their unitary invariance, that is, $\tau$ is a trace on a $B(H)$-ideal $I$ when it is a linear functional for which $\tau(UTU^*) = \tau(T)$ for all $T \in I, U \in \mathcal U(H)$. 
And essential for this is that $Ad_ U$ preserves $I$, that is, for every $X \in I$ and $U \in \mathcal U(H), Ad_U(X) := UXU^* \in I$. 
But for $J$-ideals $\mathcal I$, $Ad_U$ may not preserve $\mathcal I$ (Example \ref{E:1} above). 
However some adjustments can be made to preserve much of the trace notion. 

\bD{D:1}
For a $J$-ideal $\mathcal I$ and the subgroup of unitary operators $\mathcal U_{J}(H):= \{{\bf1}+A \in \mathcal U(H) |\, A \in J\}$, a linear functional
$$\tau : \mathcal I \rightarrow \mathbb C$$  is called a subideal-trace if $\tau(X) = \tau(UXU^*)$ for every $X \in \mathcal I$ and $U \in \mathcal U_{J}(H)$. 
In other words, $\tau$ is called a subideal-trace if $\tau$ is $\text{Ad}_{\mathcal U_{J}(H)}$-invariant, that is, if $\tau(X) = \tau(Ad_U(X))$ for $U \in \mathcal U_{J}(H)$ and $X \in \mathcal I$. 
\eD

\bR{R:00}
In particular, if $J = B(H)$ (so $\mathcal U_{B(H)}(H) = \mathcal U(H)$), then $\mathcal I$ is a $B(H)$-ideal and hence $Ad_{U}$  preserves $\mathcal I$ for $U \in \mathcal U(H)$ and Definition \ref{D:1} becomes the standard definition of a trace on a $B(H)$-ideal.
\eR

\bE{E:5.1}
\textbf{(A simple example of a subideal-trace)}\\
Consider $(S)_J$, a principal linear $J$-ideal  generated by $S \in J$ that is not a $B(H)$-ideal, and recall Proposition \ref{P:11} on the structure of its elements.
Define the map
$\tau:(S)_J \rightarrow \mathbb C $ as 
$$\tau(\alpha S + AS + SB + \displaystyle{\sum_{k=1}^{m}}A_kSB_k) := \alpha,$$ 
where $A, B, A_k, B_k \in J, \alpha \in \mathbb C, m \in \mathbb N$.
By our methods developed earlier, it is elementary to show that $\tau$ is a well-defined linear functional on $(S)_J$ when $(S)_J$ is not a $B(H)$-ideal.
Indeed, if $\alpha S + X = \beta S + Y$ for $X, Y \in SJ + JS + J(S)J$, then $(\alpha - \beta)S \in SJ + JS + J(S)J$. 
Since $(S)_J$ is not a $B(H)$-ideal, $\alpha = \beta$ (otherwise $S \in J(S)$ which by Theorem \ref{T:3.2.2} implies $(S)_J$ is a $B(H)$-ideal).
Therefore $\tau(\alpha S + X) = \tau(\beta S + Y)$, hence $\tau$ is a well-defined map.
It is elementary to show that $\tau$ is a linear map.
And since $$({\bf 1}+ A)(\alpha S + AS + SB + \displaystyle{\sum_{k=1}^{m}}A_kSB_k)({\bf 1}+ A^*) = \alpha S + X ~\text{for}~ X \in SJ + JS + J(S)J,$$ it follows that this $\tau$ is $\text{Ad}_{\mathcal U_{J}(H)}$-invariant. 
Hence $\tau$ is a subideal-trace on $(S)_J$.
\eE

The commutator space of a $B(H)$-ideal $I$, $[I, B(H)]$, is the linear span of single commutators $[A, B]$ for $A \in I, B \in B(H)$.
Since $UXU^* - X = [UX, U^*] \in [I, B(H)]$ for every $X\in I$ and every unitary operator $U \in \mathcal U(H)$, and since unitary operators span $B(H)$, unitarily invariant linear functionals on $I$ are precisely the linear functionals on $I$ that vanish on the commutator space $[I, B(H)]$ \cite[Section 2]{KW10}.

Because every operator is the linear combination of four unitary operators, the well-known commutator space $[I, B(H)]$ is also the linear span of the single commutators $[A, U]$ for $A \in I, U \in \mathcal U(H)$.
That is, $[I, \mathcal U(H)] = [I, B(H)]$. 
Observing that $\mathcal U_{B(H)}(H) = \mathcal U(H)$, we make the following analog.

\bD{D:2}
The  $\mathcal U_{J}(H)$-commutator space of $J$-ideal $\mathcal I$ is defined as $$[\mathcal I, \mathcal U_{J}(H)] := \text{linear span}\{[X,U]\,|\, X \in \mathcal I, U \in \mathcal U_{J}(H)\}$$
\eD

Notice that if $\mathcal I$ is a $B(H)$-ideal, then the  $\mathcal U_{J}(H)$-commutator space of $\mathcal I$ is precisely $[\mathcal I, B(H)]$, the commutator space of $\mathcal I$.

In the following proposition we obtain a necessary and sufficient condition for a linear functional on a subideal to be a subideal-trace. This is an analog of the trace case just described.

\bP{P:1}
For a $J$-ideal $\mathcal I$, a linear functional $\tau : \mathcal I \rightarrow \mathbb C$ is a subideal-trace if and only if $\tau$ vanishes on the $\mathcal U_{J}(H)$-commutator space of $\mathcal I$, that is, $\tau$ vanishes on $[\mathcal I, \mathcal U_{J}(H)]$.
\eP

\bp
Suppose $\tau$ is a subideal-trace.
It suffices to show that $\tau$ vanishes on single commutators $[X, U]$ for $X \in \mathcal I$ and $U \in \mathcal U_{J}(H)$.
For $X \in \mathcal I$ and ${\bf1} + B \in \mathcal U(H)$ where $B \in J$, $X({\bf1} + B)= X + XB \in \mathcal I$ .
Since $\tau$ is $Ad_{\mathcal U_{J}(H)}$-invariant, 
\begin{align*}
\tau(X({\bf1} + B)) &= \tau(({\bf1} + B)X({\bf1} + B)({\bf1} + B^*))\\
&= \tau(({\bf1} + B)X))\\
\tau(X({\bf1} + B)-({\bf1} + B)X) &=0\\
\tau([X, ({\bf1} + B)])&=0
\end{align*}
Therefore  $\tau([X, U]) = 0$ for every $U \in \mathcal U_{J}(H)$.

Next we prove the reverse implication, that is, if $\tau$ vanishes on the $\mathcal U_{J}(H)$-commutator space of $\mathcal I$, $[\mathcal I, \mathcal U_{J}(H)]$, then $\tau$ is a subideal-trace.
That is, for $U \in \mathcal U_{J}(H)$, $\tau(X) = \tau(UXU^*)$.

Since $\tau$ vanishes on $[\mathcal I, \mathcal U_{J}(H)]$, in particular, $\tau([X, ({\bf1} + B)]) =0$ implying $\tau(BX) = \tau(XB)$ for all $X \in \mathcal I$ and $({\bf1} + B) \in \mathcal U_J(H)$. 
Since $U = {\bf1} + B$ is a unitary operator, $({\bf1} + B)({\bf1} + B^*) = {\bf1}$ hence $B + B^* + BB^* = 0$.
\begin{align*}
\tau(({\bf1} + B^*)X({\bf1} + B))- \tau(X)&= \tau((X + B^*X)({\bf1} + B)-X)\\
&= \tau(X+XB+B^*X+ B^*XB - X)\\
&= \tau(XB+B^*X+ B^*XB)\\
&= \tau(XB) + \tau(B^*X) + \tau(B^*XB)\\
&= \tau(BX) + \tau(B^*X) + \tau(BB^*X) \qquad(\text{since}~ B^*X \in \mathcal I)\\
&= \tau(BX + B^*X + BB^*X)\\
&= \tau((B+B^*+BB^*)X) = \tau(0) = 0
\end{align*}

\noindent Therefore linear functional $\tau$ is $Ad_{\mathcal U_{J}(H)}$-invariant, and hence by Definition  \ref{D:1}, $\tau$ is a subideal-trace on $\mathcal I$.
\ep

\bC{C:1}
The set of all subideal-traces on a $J$-ideal $\mathcal I$ can be identified with the elements of the linear dual of the quotient space $\frac{\mathcal I}{[\mathcal I,  \mathcal U_J(H)]}$.

Indeed, for a given subideal-trace $\tau$ on a subideal $\mathcal I$, define a functional $f_{\tau} :  \frac{\mathcal I}{[\mathcal I,  \mathcal U_J(H)]} \rightarrow  \mathbb C$ as $f_{\tau}([X]):= \tau(X)$ where $[X]$ is the coset of the element $X \in \mathcal I$. 
Since $[X] = [Y]$ implies $X-Y \in [\mathcal I,  \mathcal U_J(H)]$ and $\tau$ a subideal-trace, $\tau(X-Y) = 0$ which implies that $f_{\tau}$ is a well-defined linear functional on the quotient space.
On the other hand, given a linear functional $f$ on the quotient space $\frac{\mathcal I}{[\mathcal I,  \mathcal U_J(H)]}$, define a function $\tau : \mathcal I \rightarrow \mathbb C$ as $\tau(X) := f([X])$. 
Since $f$ is a linear functional, $\tau$ is also a linear functional. 
And for every element $Y \in [\mathcal I,  \mathcal U_J(H)], f([Y]) = 0$ implying $\tau(Y) = 0$. 
Hence $\tau$ vanishes on $[\mathcal I,  \mathcal U_J(H)]$.
Therefore by Proposition \ref{P:1}, $\tau$ is a subideal-trace on $\mathcal I$.
\eC

\bR{R:111} A subideal $\mathcal I$ may be invariant under a larger class than $U \in \mathcal U_{J}(H)$ but not invariant under the full group of unitary operators  $\mathcal U(H)$. 
For instance, $U = \lambda ({\bf1} + B)$ for $|\lambda| = 1$ and $({\bf1} + B) \in \mathcal I$.
But there may be more less obvious unitary operators under which $\mathcal I$ is invariant (Example \ref{example} below).
This leads us to suggest the following alternative definition of a subideal-trace (Definition \ref{Def} below). 
However we will not explore it further here.
\eR
 
\bE{example}\quad \\
(A $K(H)$-ideal invariant under a larger class of unitaries, but not invariant under the full group $\mathcal U(H)$)

Using the principal $K(H)$-ideal $(S)_{K(H)}$ and the unitary map $\phi$ of Example \ref{E:1}, 
the unitary operator $U := \phi^{-1}({\bf1}\oplus({-\bf1}))\phi \in \mathcal U(H)\setminus \mathcal U_{K(H)}(H)$. 
That $U \notin \mathcal U_{K(H)}(H)$ is a simple computation.
Then $(S)_{K(H)}$ is invariant under $Ad_U$ because $USU^* = S$ (an easy verification combining the definition of $U$ here with the definition of $S$ in Example \ref{E:1}), but $(S)_{K(H)}$ is not invariant under $Ad_U$ for $U \in \mathcal U(H)$ which again follows from Example \ref{E:1}.  
\eE

\bD{Def}
For a $J$-ideal $\mathcal I$ and $\mathcal U^{\mathcal I}(H):=\{U \in \mathcal U(H)|\, UXU^* \in \mathcal I\,\, \text{for}\,\, X \in \mathcal I\}$, a linear functional
$$\tau : \mathcal I \rightarrow \mathbb C$$ is called a $\mathcal U^{\mathcal I}(H)$-subideal-trace if $\tau(X) = \tau(UXU^*)$ for every $X \in \mathcal I$ and $U \in \mathcal U^{\mathcal I}(H)$, that is, $\tau$ is $\text{Ad}_{\mathcal U^{\mathcal I}(H)}$-invariant. 
\eD 

The following inclusion holds for a subideal $\mathcal I$:
\begin{center}
$Ad_{\mathcal U^{\mathcal I}(H)}$-invariant subideal-traces of $\mathcal I$ $\subset$ $Ad_{\mathcal U_{J}(H)}$-invariant subideal-traces of $\mathcal I$ 
\end{center}

The next natural question is whether or not these inclusions are proper. 
In particular, do Definition \ref{D:1} and Definition \ref{Def} define different classes of functionals on a subideal that is not a $B(H)$-ideal? 
When $\mathcal I$ is a $B(H)$-ideal, Remark \ref{R:00} tells us that they are the same class.

\end{document}